\documentclass[12pt]{article}

\usepackage{amsmath}
\usepackage{amssymb}
\usepackage{amsfonts}
\usepackage{lscape}

\begin{document}

\fontsize{12}{6mm}\selectfont
\setlength{\baselineskip}{2em}

$~$\\[.35in]
\newcommand{\dss}{\displaystyle}
\newcommand{\raro}{\rightarrow}
\newcommand{\be}{\begin{equation}}

\def\sech{\mbox{\rm sech}}
\def\sn{\mbox{\rm sn}}
\def\dn{\mbox{\rm dn}}
\thispagestyle{empty}

\begin{center}
{\Large\bf Connections of Zero Curvature } \\ [2mm]   
{\Large\bf and Applications to Nonlinear }  \\   [2mm]
{\Large\bf Partial Differential Equations}  \\    [2mm]
\end{center}

\vspace{1cm}
\begin{center}
\noindent
{\bf Paul Bracken}            \\            
{\bf Department of Mathematics,}   \\
{\bf University of Texas,}   \\
{\bf Edinburg, TX  }     \\
{78540-2999}
\end{center}

\vspace{2cm}
\begin{abstract}
A general formulation of zero curvature connections
in a principle bundle is presented and some applications are
discussed. It is proved that a related connection based on a
prolongation in an associated
bundle remains zero curvature as well. It is also shown that 
the connection coefficients can be defined so that the
partial differential equation to be studied appears as 
the curvature term
in the structure equations. It is discussed how
Lax pairs and B\"acklund tranformations
can be formulated for such equations that occur
as zero curvature terms.
\end{abstract}

\vspace{1cm}
MSCs: 53Z05, 57R99, 55R10, 53B50

\vspace{4mm}
\newpage
\section{Introduction}
\numberwithin{equation}{section}

Connections which determine representations of zero curvature 
have turned out to be a very useful and innovative
approach for studying nonlinear partial differential equations.
These connection forms have the capacity 
to produce results which can be used
to obtain Lax pairs as well as B\"acklund
transformations in a very direct way provided 
information concerning the structural differential 
forms of special fiber bundles
can be specified. These types of connection have a special 
property in that the curvature tensor of such a connection
contains a subtensor which is directly proportional to 
a partial differential equation which is of interest.
For the case in which the connection tensor
with these components vanishes, as on the corresponding lifts
of solutions of a given nonlinear equation, it is said the
connection determines a representation of zero curvature.

The main ideas which have led to these developments began
several decades ago and can be traced to the work of
people such as Estabrook and Wahlquist {\bf [1-4]} and by R. Hermann {\bf [5]}
as well. Hermann first introduced at one point a particular
connection of basically this type. He proposed early on to
interpret the B\"acklund transformation as a connection
similar in a certain sense to the connection which
defines a representation of zero curvature. He first 
introduced the concept of a B\"acklund connection which is defined
by the way the connection form is specified. Hermann then
formulates B\"acklund's problem as that of finding
a section in a bundle space on whose
pull-back the B\"acklund connection is plane. He has
presented the basic idea in {\bf [6]}, and an 
introductory outline can be given based on that.

Let $M$ be a manifold and consider two sorts of object
on $M$. First $I$ will be a differential ideal of
differential forms on $M$, and $R$ a Pfaffian system
or submodule of the set of differential one-forms on $M$. 
Thus, $F^* (M)$ denotes the exterior algebra of
differential forms on $M$, and $R$ is called a
prolongation of $I$ if the following condition is
satisfied
\be
d R \subset F^{*} (M) \wedge R + I.
\label{eqI1}
\end{equation}
In the initial approach taken by Estabrook and Wahlquist,
they primarily start off with $I$ and then search for $R$.
If $I=0$, then \eqref{eqI1} expresses the fact that $R$
is completely integrable. The Frobenius complete integrability
theorem {\bf [7]} then asserts that there are, locally, one-forms
$\omega_1, \cdots, \omega_n \in R$ forming a basis and
such that $d \omega_1 = \cdots = d \omega_n =0$.
Second, if $R$ is generated by a single element, $\omega$,
such that $d \omega \in I$, then $\omega$ is a conservation
law for $I$. Studying the relation \eqref{eqI1} in more
advanced ways and further generalizations has led to an
entire geometric approach to the classic AKNS program
{\bf [8-9]}, and the study of the geometric properties of
non-linear partial differential equations and their
associated solutions. There has been much interest
in this approach {\bf [10-13]}, and has led to many
insights between integrable evolution equations and
pseudo-spherical surfaces as well {\bf [14-16]}.

The objective of this work is to go beyond this more
primitive formulation which has just been described by
starting with a jet-bundle $J^r E$ of $r$-jets over
a lower dimensional bundle $E$ {\bf [17]}. For purposes here, 
$r$ is usually two or three when second or third order
equations are involved, however, a formulation which
doesn't specify $r$ at first will be given.
Structure equations are established for the systems
of forms on these bundles. A very novel approach to
the formulation of zero curvature connections is
presented in detail. Several theorems and different
proofs of these are presented as well which establish a general 
theory of the subject from a specific abstract viewpoint. 
It is shown how the choice of particular connection
coefficients can lead to an expression for the curvature,
and an expression for the curvature tensor under the assumed
form of the coefficients is found and satisfies a particular relation.
It is also shown how prolongations of the connections can be
generated, and the resulting connections remain zero curvature.
Out of this comes a method for writing Lax pairs and 
B\"acklund transformations {\bf [18]} for the equations involved.
In fact, one of the remarkable
features of these differential systems is that
once they have been specified, they can 
be used to yield Lax pairs very easily as well as B\"acklund
transformations for the equations which appear as
the zero curvature terms in the 
structure equations.
It is explained in detail how these can be constructed. The
difficult part as far as applications are concerned
is to be able to write down the specific system
of connection one-forms to initialize the process. These same forms
contain the relevant information for producing these 
additional structures. Finally, it will
be shown how the formalism can be applied in
practice to obtain B\"acklund transformations between the
Liouville equation and the wave equation.
Differential systems which are the zero
curvature representations for these two different nonlinear
equations will be written down. They will be shown
to have the right zero curvature structure and moreover
how information  
from these differential forms needed to
write down Lax pairs and B\"acklund transformations can be extracted.

\section{Geometrical Setting}

\subsection{Framework}

The main purpose in formulating connections
which define representations of zero curvature
is to study nonlinear partial differential equations
in a systematic way. By this it is intended
that useful structures relevant to the study of
these equations,
such as Lax pairs and B\"acklund transformations,
can be produced. For definiteness, a general
third order equation is of the form
\be
F ( x^i, u, u_j, u_{jk}, u_{jkl}) =0.
\label{eqII1}
\end{equation}
By enlarging the manifold which supports \eqref{eqII1},
equations of this type can be written in a more general
form as
\be
F ( x^i, u, \lambda_j, \lambda_{jk}, \lambda_{jkl})=0,
\label{eqII2}
\end{equation}
This notation is common and can be found in {\bf [19-20]}.
The $\{ x^i, u\}$ are adapted local coordinates in
the $(n+1)$-dimensional bundle $E$ over the
$n$-dimensional base $M$, whose local coordinates
are given by $\{ x^i \}$ where $i,j,k= 1, \cdots, n$.
This larger manifold called $J^r E$ over which \eqref{eqII2}
is defined is called the space of
holonomic $r$-jets of the local sections of the
manifold $E$. It carries the system of coordinates
$\{ x^i$, $u$, $\lambda_{j_1, \cdots, j_k} \}$ with $k=1, \cdots, r$.
Thus, there exist the following
inclusions, $M \subset E \subset J^r E$.
Let $\omega^i$, $\omega^{n+1}$, $\omega_j^i$,
$\omega^{n+1}_j$, $\omega^{n+1}_{n+1}$, $\omega^i_{jk},
\cdots$ be a sequence of structural forms of
the holonomic frames of the manifold $E$, symmetric
in the subscripts. The forms $\omega^i$, $\omega^{n+1}$,
$\omega^{n+1}_{i_1, \cdots, i_k}$, for $k=1, \cdots, r$,
are referred to as principal forms in the bundle
of holonomic $r$-jets, $J^r E$ {\bf [21]}. These forms will
satisfy systems of structural equations which 
have the form,
\be
\begin{array}{c}
d \omega^i = \omega^k \wedge \omega_j^i,  \\
  \\
d \omega^{n+1} = \omega^j \wedge \omega^{n+1}_j 
+ \omega^{n+1} \wedge \omega_{n+1}^{n+1},  \\
\end{array}
\label{eqII3}
\end{equation}
as well as equations which arise in the process
of regular prolongation of these by means of Cartan's
lemma.
That is to say, taking the exterior derivative of
the first equation in \eqref{eqII3} gives
$$
0 = d^2 \omega^i = d \omega^k \wedge \omega_k^i
- \omega^k \wedge d \omega_k^i
= \omega^s \wedge ( \omega_s^k \wedge \omega_k^i
- d \omega_s^i).
$$
By the generalized Cartan lemma, the coefficients
in the brackets can be expanded in terms of the
forms $\omega^i$
$$
d \omega_s^i - \omega_s^k \wedge \omega_k^i
= \omega^k \wedge \omega^i_{sk}.
$$
This can be differentiated in turn and when the
process is repeated, a tower of forms can be 
constructed {\bf [22]}.

It is important in the course of this work to be able
to evaluate appropriate sections in these bundles,
and it is carried out in the following way.
For any section $\Sigma \subset E$ which is defined
by the equation $u = u (x^1, \cdots, x^n)$, 
sections in $\Sigma^r \subset J^r E$ are defined by 
the equations
\be
u= u(x^1, \cdots, x^n),
\qquad
\lambda_{i_1, \cdots, i_k} = u_{i_1, \cdots, i_k},
\qquad
k=1, \cdots, r.
\label{eqII4}
\end{equation}
The subscripts $i+1, \cdots, i_k$ on the function $u$
now denote partial derivatives. Consequently, under
this process, the equation \eqref{eqII2} is mapped onto
\eqref{eqII1}, the equation of interest.
If contact forms are chosen as principal forms on
the manifold $J^r E$, then the pull-backs are integral
manifolds of the system of Pfaffian equations
\be
\omega^{n+1} = \omega^{n+1}_i = \cdots = \omega^{n+1}_{i_1 \cdots i_k} 
=0.
\label{eqII5}
\end{equation}

\subsection{Principle Bundle}

To begin with, based on this sequence of manifolds, consider 
the principle bundle $P (J^r E, G)$ over $J^r E$ along with the
$g$ parameter structure group $G$. 
Let $P (J^r E, G)$ have structural forms $\omega^A$,
($A, B =1, \cdots, g$) which satisfy structure equations of
the form
\be
d \omega^A = \frac{1}{2} C_{BC}^A \, \omega^B \wedge \omega^C
+ \omega^{\delta} \wedge \omega_{\delta}^A.
\label{eqII6}
\end{equation}
In \eqref{eqII6}, the $C_{BC}^A$ are the structure constants
pertaining to the Lie group $G$. They are skew-symmetric with
respect to the lower indices and satisfy the Jacobi identity
\be
C^A_{BK} C_{LM}^B + C_{BL}^A C_{MK}^B + C_{BM}^A C_{KL}^B =0.
\label{eqII7}
\end{equation}
The forms $\omega^{\delta}$ will be principle forms of the
base $J^r E$, and will be completely integrable.
Thus, their differentials satisfy structure equations
of the form
\be
d \omega^{\delta} = \omega^{\mu} \wedge \omega_{\mu}^{\delta}.
\label{eqII8}
\end{equation}

\section{General Zero-Curvature Formulation}

To show exactly how zero curvature representations can 
be developed from a rigorous point of view, a connection
in the principle bundle $P (J^r E, G)$ has to be defined {\bf [19-20]}.
One way of doing this is to specify the object of
connection. This is made precise in the following theorem.

{\bf Theorem 3.1} A connection in the principle bundle
$P (J^r E, G)$ can be given by the field of a connection
object on $J^r E$ which has components $\Gamma^A_{\epsilon}$
that satisfy the system of differential equations
\be
d \Gamma^A_{\epsilon} + C_{BC}^A \Gamma^B_{\epsilon} \omega^C
- \Gamma^A_{\delta} \omega^{\delta}_{\epsilon} - \omega_{\epsilon}^A
= \Gamma_{\epsilon \delta}^A \omega^{\delta},
\label{eqIII1}
\end{equation}
The forms $\omega_{\epsilon}^{\delta}$ are determined from
\eqref{eqII8}. The associated connection forms
\be
\tilde{\omega}^A = \omega^A + \Gamma^A_{\epsilon} \, \omega^{\epsilon}
\label{eqIII2}
\end{equation}
satisfy the structure equations
\be
d \tilde{\omega}^A = \frac{1}{2} C_{BC}^A \tilde{\omega}^B \wedge 
\tilde{\omega}^C + \Omega^A.
\label{eqIII3}
\end{equation}
The $\Omega^A$ in \eqref{eqIII3} are curvature forms given by
\be
\Omega^A = R^A_{\epsilon \delta} \omega^{\epsilon} \wedge
\omega^{\delta}.
\label{eqIII4}
\end{equation}

{\bf Proof:} Differentiating the connection forms in \eqref{eqIII2}
and requiring the exterior derivative be consistent with \eqref{eqIII3},
yields
$$
d \omega^A + d ( \Gamma^A_{\delta} \omega^{\delta})
= \frac{1}{2} C_{BC}^A ( \omega^B + \Gamma^B_{\epsilon} \omega^{\epsilon})
\wedge ( \omega^C + \Gamma^C_{\delta} \omega^{\delta}) + \Omega^A.
$$
Expanding this out, the following expression results,
$$
d \omega^A + d \Gamma^A_{\delta} \wedge \omega^{\delta} + \Gamma_{\delta}^A
d \omega^{\delta} = \frac{1}{2} C_{BC}^A \omega^B \wedge \omega^C +
\frac{1}{2} C_{BC}^A \omega^B \wedge \Gamma^C_{\delta} \omega^{\delta}
+ \frac{1}{2} C_{BC}^A \Gamma^B_{\epsilon} \omega^{\epsilon} \wedge
\omega^C + \frac{1}{2} C_{BC}^A \Gamma_{\epsilon}^B \Gamma^C_{\delta}
\omega^{\epsilon} \wedge \omega^{\delta} + \Omega^A.
$$
Substituting \eqref{eqII8} and \eqref{eqIII1} into this, we obtain,
$$
d \omega^A - \frac{1}{2} C_{BC}^A \omega^B \wedge \omega^C - \omega^{\delta}
\wedge \omega^A_{\delta} + ( - C_{BC}^A \Gamma^B_{\delta} \omega^C
+ \Gamma^A_{\sigma} \omega^{\sigma}_{\delta} + \omega^A_{\delta}
+ \Gamma_{\delta \sigma}^A \omega^{\sigma}) \wedge \omega^{\delta}
+ \Gamma^A_{\delta} \omega^{\epsilon} \wedge \omega^{\delta}_{\epsilon}
$$
$$
= - \omega^{\delta} \wedge \omega^A_{\delta} + \frac{1}{2} C_{BC}^A
\Gamma^C_{\delta} \omega^B \wedge \omega^{\delta} + \frac{1}{2}
C_{BC}^A \Gamma^C_{\delta} \omega^{\delta} \wedge \omega^B
+ \frac{1}{2} C_{BC}^A \Gamma^B_{\epsilon} \Gamma^C_{\delta} 
\omega^{\epsilon} \wedge \omega^{\delta} + \Omega^A.
$$
Now replace $d \omega^A$ using \eqref{eqII6} to obtain
$$
- C_{BC}^A \Gamma^B_{\delta} \omega^C \wedge \omega^{\delta}
+ \Gamma^A_{\sigma} \omega^{\sigma}_{\delta} \wedge \omega^{\delta} + \omega^A_{\delta} 
\wedge \omega^{\delta} + \Gamma^A_{\delta \sigma} \omega^{\sigma}
\wedge \omega^{\delta} + \Gamma^A_{\delta} \omega^{\epsilon} 
\wedge \omega_{\epsilon}^{\delta}
$$
$$
=- \omega^{\delta} \wedge \omega^A_{\delta} + C_{BC}^A 
\Gamma^C_{\delta} \omega^B \wedge \omega^{\delta} + \frac{1}{2}
C_{BC}^A \Gamma^B_{\epsilon} \Gamma^C_{\delta} \omega^{\epsilon}
\wedge \omega^{\delta} + \Omega^A.
$$
The fact that the $C_{BC}^A$ are antisymmetric in the lower
indices simplifies this result to the form,
$$
\Omega^A = \Gamma^A_{\delta \sigma} \omega^{\sigma} \wedge
\omega^{\delta} - \frac{1}{2} C_{BC}^A \Gamma^B_{\epsilon} 
\Gamma^C_{\delta} \, \omega^{\epsilon} \wedge \omega^{\delta}.
$$
Factoring the one-forms in the first part of $\Omega^A$,
it is found that
\be
\Omega^A =- \frac{1}{2} ( \Gamma_{\epsilon \delta}^A 
- \Gamma_{\delta \epsilon}^A + C_{BC}^A \Gamma^B_{\epsilon}
\Gamma_{\delta}^C) \, \omega^{\epsilon} \wedge \omega^{\delta}.
\label{eqIII5}
\end{equation}
This gives $\Omega^A$ explicitly and finishes the proof.

The coefficients of $\Omega^A$ in \eqref{eqIII5}
give the components of $R^A_{\epsilon \delta}$ and the
theorem allows us to identify the components of the
curvature tensor as
\be
R^A_{\epsilon \delta} =- \frac{1}{2} ( \Gamma^A_{\epsilon \delta}
- \Gamma^A_{\delta \epsilon} + C_{BC}^A \Gamma^B_{\epsilon}
\Gamma^C_{\delta} ).
\label{eqIII6}
\end{equation}

{\bf Theorem 3.2} The curvature tensor satisfies the following
relation
\be
d R_{\lambda \mu}^A + R_{\lambda \mu}^B C_{BC}^A \omega^C
- R_{\sigma \mu}^A \omega_{\lambda }^{\sigma}
- R_{\lambda \sigma}^A \omega_{\mu}^{\sigma} =0,
\qquad \mod \, \omega^{\Delta},
\label{eqIII7}
\end{equation}
where $\omega^{\Delta}$ are principle forms of the jet
manifold.

{\bf Proof:} Differentiating both sides of \eqref{eqIII3}
exteriorly, it is found that
$$
0= \frac{1}{2} C_{BC}^A \, d \tilde{\omega}^B \wedge \tilde{\omega}^C
- \frac{1}{2} C_{BC}^A \tilde{\omega}^B \wedge d \tilde{\omega}^C
+ d R_{\lambda \mu}^A \wedge \omega^{\lambda} \wedge \omega^{\mu}
+ R_{\lambda \mu}^A d \omega^{\lambda} \wedge \omega^{\mu}
- R_{\lambda \mu}^A \, \omega^{\lambda} \wedge d \omega^{\mu}
$$
$$
= C_{BC}^A ( \frac{1}{2} C_{DQ}^B \tilde{\omega}^D \wedge \tilde{\omega}^Q
+ R_{\lambda \mu}^B \omega^{\lambda} \wedge \omega^{\mu}) \wedge \tilde{\omega}^C
+ d R_{\lambda \mu}^A \wedge \omega^{\lambda} \wedge \omega^{\mu}
+ R_{\lambda \mu}^A \omega^{\sigma} \wedge \omega_{\sigma}^{\lambda} 
\wedge \omega^{\mu} - R_{\lambda \mu}^A \, \omega^{\lambda} \wedge
\omega^{\sigma} \wedge \omega_{\sigma}^{\mu}
$$
$$
= \frac{1}{2} C_{TC}^A C_{DB}^T \tilde{\omega}^D \wedge \tilde{\omega}^B
\wedge \tilde{\omega}^C + C_{BC}^A R_{\lambda \mu}^B \tilde{\omega}^C
\wedge \omega^{\lambda} \wedge \omega^{\mu} + d R_{\lambda \mu}^A \wedge \omega^{\lambda}
\wedge \omega^{\mu} - R_{\lambda \mu}^A \omega_{\sigma}^{\lambda} \wedge
\omega^{\sigma} \wedge \omega^{\mu} - R_{\lambda \mu}^A \, \omega_{\sigma}^{\mu}
\wedge \omega^{\lambda} \wedge \omega^{\sigma}.
$$
Invoking the Jacobi identity \eqref{eqII7}, this result reduces to the following form
$$
( d R_{\lambda \mu}^A + R_{\lambda \mu}^B C_{BC}^A \tilde{\omega}^C
- R_{\sigma \mu}^A \omega_{\lambda}^{ \sigma} - R_{\lambda \sigma}^A
\omega_{\mu}^{\sigma} ) \wedge \omega^{\lambda} \wedge \omega^{\mu} =0.
$$
This implies that the coefficient of $\omega^{\lambda} \wedge \omega^{\mu}$
is zero $\mod \omega^{\Delta}$, the principle forms of the jet manifold, so that
$\tilde{\omega}^C = \omega^C$. The result in \eqref{eqIII7} then follows.

Thus, the curvature tensor components include, in particular,
the components $R_{kl}^A$. As a consequence of these theorems, the
following result is very important as far as the application of the
zero-curvature idea to specific nonlinear differential equations
is concerned.

{\bf Theorem 3.3} For the connection given in the principle bundle
$P (J^r E, G)$ to define the representation of zero curvature
which corresponds to an equation $F (x^i, u, \lambda_j, \lambda_{jk}, \cdots)=0$,
it is necessary and sufficient that the components $R_{kl}^{A}$ of the
curvature vanish on the pull-backs of the solutions to the equation.

{\bf Proof:} Since the vanishing of the forms of curvature $\Omega^A
= R_{\lambda \mu}^A \, \omega^{\lambda} \wedge \omega^{\mu}$ on the
pull-backs of solutions is invariant, it suffices to show the
statement for some special choice of the principle forms. The
statement then becomes obvious if contact forms are taken as
principle forms since, in this case, the relations $\Omega^A
= R^A_{kl} \omega^k \wedge \omega^l$ hold on the pull-back
of any section $\Sigma \subset E$.
 
In practical terms, the curvature tensor will be, or will have
a subtensor, which is proportional to the equation under
consideration, and will clearly vanish identically on solutions
of that equation. Thus, a connection is called a connection
determining a representation of zero curvature for a differential
equation if the curvature form vanishes on the solutions, or on
the corresponding lifts of solutions, and only on solutions.

\section{Prolongations on These Spaces}

An additional bundle associated with the principle
bundle $P (J^r E, E)$, which is called $F (P (J^r E, G))$,
can now be constructed. A
larger space is now being associated with $P$. The typical
fiber of this new bundle is a space $F$ which is an
$N$-dimensional space of the representation of the Lie group
$G$. The representation of the group $G$ as a group of
transformations of the space $F$ can be defined by the
specification of the system of Pfaffian equations
\be
d X^I - \xi_A^I (X) w^a=0.
\label{eqIV1}
\end{equation}
In \eqref{eqIV1}, the $w^A$ are invariant forms of the group
$G$ which satisfy the structural equations
\be
d w^A = \frac{1}{2} \, C_{BC}^A \, w^B \wedge w^C.
\label{eqIV2}
\end{equation}
Indeed, it is worth recalling that if $G$ is connected, any
diffeomorphism $f : G \raro G$ which preserves left-invariant
forms, $\theta^{\alpha}$, so that $f^* \theta^{\alpha} = \theta^{\alpha}$
is left translation. If $N$ is a smooth manifold and $w^{\alpha}$
linearly independent forms on $N$ satisfying \eqref{eqIV2}, then
for any point in $N$, there exists a neighborhood $U$ and a
diffeomorphism $f: U \raro G$ such that $\theta^{\alpha} = 
f^* (w^{\alpha})$.

The following theorem will produce a condition that, when satisfied,
will guarantee that system \eqref{eqIV1} is completely integrable.

{\bf Theorem 4.1} Pfaffian system \eqref{eqIV1} is completely
integrable provided the set of $\xi^I_A (X)$ satisfy the following 
constraint,
\be
\xi_B^K \frac{\partial \xi_C^I}{\partial X^K}
- \xi^K_C \frac{\partial \xi^I_B}{\partial X^K}
+ \xi_A^I C_{BC}^A =0.
\label{eqIV3}
\end{equation}

{\bf Proof:} Differentiate both sides of system \eqref{eqIV1} to
obtain,
$$
\frac{\partial \xi_A^I}{\partial X^K} \xi_C^K (X) \, w^C \wedge w^A
+ \frac{1}{2} \xi_B^I (X) C_{BC}^A \, w^B \wedge w^C =0.
$$
The first term in this equation can be put in the form
$$
\frac{1}{2} \{ \xi_B^K (X) \frac{\partial \xi_C^I}{\partial X^K} \,
w^B \wedge w^C + \xi_C^K (X) \frac{\partial \xi_B^I}{\partial X^K}
\, w^C \wedge w^B \} + \frac{1}{2} \xi_A^I (X) C_{BC}^A \, w^B \wedge w^C =0.
$$
Equating the coefficient of $w^B \wedge w^C$ to zero, the condition
\eqref{eqIV3} for complete integrability is obtained. These conditions
are often referred to as the Lie identities.

If there exists a connection in $P (J^r E, G)$ 
which determines a representation of zero curvature, it is remarkable that
the same property holds in the associated bundle $F(P(J^r E,G))$.
The $N$-dimensional space $F$ is coordinatized by means of
coordinates $ \{ X^i \}_1^N$ and carries a representation of
the group. Moreover, the curvature forms of $F ( P( J^r E, G))$
are defined by
\be
\theta^I = d X^I - \xi^I_A (X^1, \cdots, X^N) \omega^A,
\quad I,J,K =1, \cdots , N.
\label{eqIV4}
\end{equation}
In \eqref{eqIV4}, the $\omega^A$ are structural forms
of the principle bundle.

If a connection with the connection forms
\be
\tilde{\omega}^A = \omega^A + \Gamma^A_{\lambda} \omega^{\lambda},
\label{eqIV5}
\end{equation}
is defined in the principle bundle, then along with this
connection in the principle bundle, a connection is induced
in the associated bundle $F (P (J^r E, G))$ and it has 
connection forms
\be
\tilde{\theta}^I = d X^I - \xi^I_A (X) \tilde{\omega}^A.
\label{eqIV6}
\end{equation}

{\bf Proposition 4.1} The Pfaffian system $\tilde{\theta}^I$ 
satisfies the system of structural equations
\be
d \tilde{\theta}^I = \tilde{\theta}^K \wedge \tilde{\theta}^I_K
- \xi^I_A (X) \, R^A_{\lambda \mu} \omega^{\lambda} \wedge
\omega^{\mu}.
\label{eqIV7}
\end{equation}
The $\xi_A^I (X)$ satisfy the Lie identities \eqref{eqIV3}
and the $\tilde{\theta}_K^I$ are given by
\be
\tilde{\theta}_K^I =- \frac{\partial \xi_A^I}{\partial X^K}
\tilde{\omega}^A.
\label{eqIV8}
\end{equation}
The $R^A_{\lambda \mu}$ are the components of the curvature
tensor defined in $P (J^r E, G)$.

{\bf Proof:} Differentiating the set of forms in \eqref{eqIV6},
it is found that
$$
d \tilde{\theta}^I =- \frac{\partial \xi_A^I}{\partial X^K} d X^K
\wedge \tilde{\omega}^A - \xi_A^I (X) d \tilde{\omega}^A
$$
$$
=- \frac{\partial \xi_A^I}{\partial X^K} ( \tilde{\theta}^K
+ \xi_C^K (X) \tilde{\omega}^C ) \wedge \tilde{\omega}^A
- \xi_A^I (X) \, d \tilde{\omega}^A
$$
$$
=- \frac{\partial \xi^I_A}{\partial X^K} \tilde{\theta}^K \wedge 
\tilde{\omega}^A - \xi^K_C (X) \frac{\partial \xi_A^I}{\partial X^K}
\tilde{\omega}^C \wedge \tilde{\omega}^A - \xi_A^I (X) \, d \tilde{\omega}^A
$$
$$
= \tilde{\theta}^K \wedge (- \frac{\partial \xi^I_A}{\partial X^K}) \tilde{\omega}^A
- \xi^K_B (X) \frac{\partial \xi^I_C}{\partial X^K} \tilde{\omega}^B \wedge
\tilde{\omega}^C - \frac{1}{2} C_{BC}^A \xi^I_A (X) \tilde{\omega}^B \wedge
\tilde{\omega}^C - \xi_A^I (X) \, R^A_{\lambda \mu} 
\omega^{\lambda} \wedge \omega^{\mu}.
$$
Assuming that the Lie identities \eqref{eqIV3} hold and $\tilde{\theta}^I_K$
are defined by \eqref{eqIV8}, the desired result \eqref{eqIV7} appears directly,
$$
d \tilde{\theta}^I = \tilde{\theta}^K \wedge \tilde{\theta}^I_K - \xi^I_A (X)
\, R^A_{\lambda \mu} \omega^{\lambda} \wedge \omega^{\mu}.
$$
Therefore, if the connection defined in the principle bundle
specifies a representation of zero curvature for an equation,
then the related connection just defined in the associated bundle
generated by it will define a representation of zero curvature
as well. Its curvature tensor $\xi^I_A R^A_{\lambda \mu}$ vanishes
on sections $\Sigma \subset E$ if and only if
the sections are solutions of the equations. This has established 
the following.

{\bf Corollary 4.1} The system of forms $\tilde{\theta}^I$ defined by 
\eqref{eqIV6} is completely integrable on the pull-backs of solutions
to the associated equation and only on these solutions.

The theoretical advantage then in introducing the general formalism
is that the $R^A_{\lambda \mu}$ can be interpreted as curvature
forms with respect to this larger manifold. This also suggests
an application for these results. It is possible that a system of 
forms $\tilde{\theta}^K$ can be found such that a set of equations
of the form \eqref{eqIV7} obtain. The curvature terms may 
automatically vanish or be proportional to some nonlinear partial
differential equation of interest which vanishes on some transverse
integral manifold of solutions. Along with B\"acklund connections on 
bundles having one-dimensional fibers, B\"acklund connections on
bundles with two-dimensional fibers can be studied; for example,
on a bundle associated to a two-dimensional vector space of the
representation of the group $Sl (2)$. This connection is often
referred to as a Lax connection as it can be made to lead directly
to formulation of Lax pairs for the equation.
In this event, the specific forms can then be used to generate both
Lax pairs and B\"acklund transformations. This will be illustrated
clearly in the following general theorem below {\bf [23]}.

Hermann used a one-form with the structure \eqref{eqIV6} for the
KdV equation and realized that it could be written in a particular way
{\bf [5]}. He inferred that the Wahlquist-Estabrook prolongation
structure could be interpreted as a type of connection.
As for the form $\tilde{\theta}$, it is a form of connection 
in a bundle with a one-dimensional typical fiber associated with
the principal bundle $P (J^r E, Sl (2))$. This connection is also
a connection defining a representation of zero curvature.
Note that a one-form is a connection form in a bundle with 
a one-dimensional typical fiber associated with the principal
bundle $P ( J^r E, Sl (2))$ if and only if
it takes the form
$$
dy - \xi (y) \tilde{\theta}_0 - \xi_1^2 (y) \tilde{\theta}_1
- \xi_2^1 (y) \tilde{\theta}_2.
$$
The Lie identities satisfied by these coefficients are obtained
from the system
$$
\frac{\partial \xi_B^I}{\partial y^K} \xi_C^K 
- \frac{\partial \xi_C^I}{\partial y^K} \xi_B^K
= \xi_A^I C^A_{BC}.
$$

Consider a B\"acklund mapping in the one-dimensional case.
In this case the system of Pfaff equations that define the
B\"acklund mapping consist of a single equation
\be
dy - \xi (y) \tilde{\omega} - \xi_1^2 (y) \tilde{\omega}^1_2
- \xi^1_2 (y) \tilde{\omega}^2_1 =0.
\label{eqIV9}
\end{equation}
The Lie identities satisfied by the coefficients in \eqref{eqIV9}
are of the following form
$$
\xi \frac{\partial \xi_1^2}{\partial y} - \xi^2_1
\frac{\partial \xi}{\partial y} = \xi_1^2,
$$
\be
\xi \frac{\partial \xi^1_2}{\partial y} - \xi_2^1
\frac{\partial \xi}{\partial y} =- \xi_2^1,
\label{eqIV10}
\end{equation}
$$
\xi_1^2 \frac{\partial \xi^1_2}{\partial y}
- \xi_2^1 \frac{\partial \xi^2_1}{\partial y}
= 2 \xi.
$$

{\bf Theorem 4.2.} The Pfaff equation \eqref{eqIV9} which defines 
the B\"acklund mapping with the associated space
of the structure group $G$ of dimension one can be represented
in either of the two forms,
\be
\begin{array}{c}
d \varphi - \tilde{\omega}^2_1 - \varphi \tilde{\omega}
+ \varphi^2 \tilde{\omega}^1_2 =0,   \\
  \\
d \psi - \tilde{\omega}^1_2 - \psi \tilde{\omega}
- \psi^2 \tilde{\omega}^2_1 =0.
\end{array}
\label{eqIV11}
\end{equation}

{\bf Proof:} Take the second equation in \eqref{eqIV10} and divide 
it by $(\xi_2^1)^2$ to obtain
$$
- \frac{\xi}{( \xi^1_2)^2} \, d \xi^1_2 + \frac{d \xi}
{\xi_2^1} = \frac{dy}{\xi^1_2}.
$$
This is equivalent to
$$
d ( \frac{\xi}{\xi^1_2}) = \frac{dy}{\xi_2^1}.
$$
Define the variable $\varphi = \xi / \xi^1_2$ and use it
in this result to give,
\be
d \varphi = \frac{dy}{\xi^1_2}.
\label{eqIV12}
\end{equation}
Dividing by $( \xi^1_2)^2$, the third equation becomes
$$
- \frac{\xi_1^2}{( \xi_2^1)^2} \, d \xi^1_2 + \frac{d \xi_1^2}{\xi_2^1}
= -2 \frac{\xi}{(\xi_2^1)^2} \, dy.
$$
Consequently, using \eqref{eqIV12},
$$
d ( \frac{\xi^2_1}{\xi^1_2}) =- 2 \frac{\xi}{\xi^1_2} \frac{dy}{\xi^1_2}
=- d \varphi^2.
$$
Thus, we can identify $- \varphi^2 = \xi^2_1 / \xi_2^1$. Since the
form \eqref{eqIV9} can be written in the following way,
\be
\frac{dy}{\xi_2^1} - \tilde{\omega}^2_1 - \frac{\xi (y)}{\xi_2^1 (y)}
\tilde{\omega} - \frac{\xi^2_1 (y)}{\xi^1_2 (y)} \tilde{\omega}^1_2 =0,
\label{eqIV13}
\end{equation}
the required first equation in \eqref{eqIV11}
follows by substituting these results for
$\varphi$ and $\varphi^2$ into \eqref{eqIV13}. 
The second equation in \eqref{eqIV11} follows in a similar fashion.

An example which shows how the results in these last two sections
can be combined and made into something useful will be presented.
Here $M$ will be the two-dimensional base manifold which is
coordinatized by the coordinates $(x^1, x^2 ) = (x, t)$.
Now consider the following application which starts with
Theorem 3.1. A system of structural forms $\tilde{\omega}^A$
is required to satisfy the structure equations \eqref{eqIII3}
expressed as
\be
d \tilde{\omega}^1 =2 \tilde{\omega}^2 \wedge \tilde{\omega}^3
+ R_{12} \, d x^1 \wedge d x^2,
\quad
d \tilde{\omega}^2 = \tilde{\omega}^1 \wedge \tilde{\omega}^2
+ R^2_{112} \, dx^1 \wedge d x^2,
\quad
d \tilde{\omega}^3 = \tilde{\omega}^3 \wedge \tilde{\omega}^1
+ R^1_{212} \, dx^1 \wedge dx^2.
\label{eqIV14}
\end{equation}
The last terms in these are the curvature terms which are
required to vanish when they are considered on the lifting
of a section. This will result in producing a particular 
equation in the end. In the notation of \eqref{eqII2}, take
for the forms $\tilde{\omega}^A$
\be
\tilde{\omega}^1 = 2 \lambda_1 \, d x^2,
\quad
\tilde{\omega}^2 = \frac{1}{2} \lambda_1 \, d x^1 +
( u \lambda_1 - \lambda_{11} ) \, d x^2,
\quad
\tilde{\omega}^3 = d x^1 + 2 u \, d x^2.
\label{eqIV15}
\end{equation}
It is easily verified that these forms satisfy system 
\eqref{eqIV14}. The curvature term in the first and third
is zero. The second is satisfied provided that 
considered on the lifting of a section 
in which the notation reverts to that of \eqref{eqII1},
$u$ satisfies the following Burgers-type
equation $-\frac{1}{2} u_{12}
+ \frac{1}{2} (u^2)_{11} - u_{111} + (u_1)^2 =0$.
Replacing $(x^1, x^2)=(x,t)$ in this, the following
form for the equation is obtained,
\be
u_{xt} = (u^2)_{xx} - 2 u_{xxx} + 2 (u_x)^2.
\label{eqIV16}
\end{equation}
Following along the lines of Theorem 4.2, there
should be a B\"acklund transformation of the form
$dy + \tilde{\omega}^2 - y \tilde{\omega}^1
- y^2 \tilde{\omega}^3 =0$. Substituting the forms
\eqref{eqIV15} into this, the following
differential system is obtained
\be
y_x = - \frac{1}{2} u_x + y^2,
\qquad
y_t = u_{xx} - u u_x + 2 u_x y + 2 u y^2.
\label{eqIV17}
\end{equation}
Evaluating the derivatives $y_{xt}$ and $y_{tx}$,
and subtracting, all higher order terms 
in the expression above $y^0$
are found to cancel. Only the $y^0$ term remains
and it is precisely the equation \eqref{eqIV16}.

Another approach to Lax and B\"acklund systems 
will be presented in the next section.

\section{Lax and B\"acklund Systems}

Perhaps the most interesting aspect of the theoretical development
presented so far is that there exists a clear relationship between
connections which define a representation of zero curvature and
specific Lax and B\"acklund systems for the equation.
Let the group be $G = Gl (2)$, so that $r$ is selected 
to suit the system under consideration.
In fact, for the example given here, we take $r=2$, and the following 
theorem holds. 

{\bf Theorem 5.1}  Given a connection in $P (J^r E, G)$, where
$G= Gl(2)$ or a subgroup, which defines a representation of 
zero curvature corresponding to an equation of the form \eqref{eqII1},
a Lax system exists which can be defined in terms of the
connection coefficients.

{\bf Proof:} Let
\be
\tilde{\omega}^i_j = \omega_j^i + \Gamma^i_{j \lambda} \omega^{\lambda}
\label{eqV1}
\end{equation}
be connection forms in the principle bundle $P (J^r E, Gl (2))$ which
define the representation of zero curvature for the $\tilde{\omega}^A$.
This connection which is defined in the principle bundle generates a
connection in the associated bundle whose typical fiber is a two-dimensional
linear space. The connection forms in the associated bundle corresponding
to the connection in $P$ can be written in the form \eqref{eqIV6}
$$
\tilde{\theta}^i = d X^i + X^j \tilde{\omega}^i_j.
$$
As for the connection in $P$, the connection in the associated bundle
is also a connection which defines a representation of zero curvature
for the equation. Consequently, the restriction of the $\tilde{\theta}^i$
to the corresponding pull-back of the section $\Sigma \subset E$ defined
by $u = u(x,y)$ is completely integrable if and only if the section
$\Sigma \subset E$ is a solution of the equation.  $\clubsuit$

In practical terms, if contact forms are taken as principle forms
then $\omega_j^i$ will be equal to zero and the forms $\tilde{\theta}^i$
take the form
\be
\tilde{\theta}^i = d X^i + X^j \Gamma^i_{j \lambda} \omega^{\lambda}.
\label{eqV2}
\end{equation}
In this case, with $(x^1, x^2)= (x,y)$, the system of equations $\tilde{\theta}^i |_{\Sigma} =0$
have the form
\be
d X^i + X^j \Gamma^i_{j1} (x,y,u,u_k,u_{kl}) \, dx
+ X^j \Gamma_{j2}(x,y,u,u_k,u_{kl}) \, dy =0.
\label{eqV3}
\end{equation}
Of course, this is equivalent to the following system of partial
differential equations
\be
X^i_x =- \Gamma^i_{j1} (x,y, u,u_k,u_{kl}) X^j,
\qquad
X^i_y =- \Gamma^i_{j2} (x,y,u,u_k,u_{kl}) X^j.
\label{eqV4}
\end{equation}
In matrix form for a two-dimensional representation of $G$,
\eqref{eqV4} can be written as
\be
\begin{pmatrix}
X^1  \\
X^2   \\
\end{pmatrix}_x =
\begin{pmatrix}
- \Gamma^1_{11}  & - \Gamma^1_{21}  \\
- \Gamma^2_{11}  & - \Gamma^2_{21}  \\
\end{pmatrix}
\begin{pmatrix}
X^1  \\
X^2  \\
\end{pmatrix}
,
\qquad
\begin{pmatrix}
X^1  \\
X^2  \\
\end{pmatrix}_y =
\begin{pmatrix}
- \Gamma^1_{12}  &  -\Gamma^1_{22} \\
- \Gamma^2_{12}  &  - \Gamma^2_{22}  \\
\end{pmatrix}
\begin{pmatrix}
X^1  \\
X^2  \\
\end{pmatrix}.
\label{eqV5}
\end{equation}

This system is completely integrable and has solutions satisfying
any initial conditions if and only if $u=u (x,y)$ is a solution of
the associated nonlinear equation. 

There are relationships between B\"acklund
transformations and the connections defining representations of
zero curvature, as Hermann pointed out {\bf [3]}. Consider restricting the
problem to investigate how to write B\"acklund transformations between 
two second order equations. Suppose $x,y,u$ and $x,y,v$ are adapted
local coordinates in bundles $E_1$ and $E_2$ respectively which
share a common base manifold $M$ with local coordinates $x$, $y$.
The variables $x,y,u, \lambda_i, \lambda_{jk}$ and $x,y,v, \mu_i,
\mu_{jk}$ are local coordinates in the bundles of second order
jets $J^2 E_1$ and $J^2 E_2$. In this case, $x,y,u, \lambda_i$
and $x,y, v, \mu_i$ are local coordinates in the corresponding
bundles of first order jets $J^1 E_1$ and $J^1 E_2$. In this
event, the equations then take the form
\be
F_1 ( x,y, u, \lambda_i, \lambda_{jk} )=0, 
\label{eqV6}
\end{equation}
and,
\be
F_2 (x, y, v, \mu_i, \mu_{jk}) =0.  \qquad  
\label{eqV7}
\end{equation}
A B\"acklund transformation between these two equations can be
defined as a system of equations
\be
\Phi ( x,y, u,v, u_i, v_j ) =0.
\label{eqV8}
\end{equation}
Equation \eqref{eqV8} will be integrable over $u$ if and only if
$v = v(x,y)$ is a solution of \eqref{eqV7} and integrable over
$v$ if and only if $u = u(x,y)$ is a solution of \eqref{eqV6}.
For any specified solution $u$ of \eqref{eqV6}, or $v$ of \eqref{eqV7},
\eqref{eqV8} makes it possible to determine a certain solution
$v$ of \eqref{eqV7}, or of \eqref{eqV6}, respectively.

It is said that a B\"acklund transformation
is established between \eqref{eqV6} and \eqref{eqV7} if
connections have been defined in the two principle bundles
$P ( J^1 E_1, G_1)$ and $P (J^1 E_2, G_2)$ which define representations
of zero curvature for each equation. In each of the manifolds
$E_1$ and $E_2$ a structure of the bundle is defined with a
one-dimensional fiber associated. In the case of $E_2$, it is
with the principle bundle $P (J^1 E_1, G_1)$ and in the case of
$E_1$ with $P ( J^1 E_2, G_2)$. Therefore, the connections which
are defined in the principle bundles and specify representations of
zero curvature generate corresponding representations of zero
curvature in the associated bundles. The forms for these
two connection forms are written ${\theta}$ and ${\vartheta}$. 

For the case in which $G_1=G_2 = Gl (2)$, the forms $\theta$
and $\vartheta$ take the form
\be
\theta = dv - \xi_j^i (v) \tilde{\omega}_i^j,
\label{eqV9}
\end{equation}
and,
\be
\vartheta = du - \eta_j^i (u) \tilde{\pi}_i^j.
\label{eqV10}
\end{equation}
The structure forms on the right of \eqref{eqV9} and 
\eqref{eqV10} are given by
\be
\tilde{\omega}_j^i = \omega_j^i + \Gamma^i_{jk} (x,y,u,\lambda_l)
\omega^k,
\quad
\tilde{\pi}^i_j = \pi_j^i + \Phi^i_{jk} (x,y,v, \mu_l) \omega^k,
\qquad i,j=1,2.
\label{eqV11}
\end{equation}
These will be connection forms in $P (J^1 E_1, Gl (2))$ and
$P ( J^1 E_2, Gl(2))$, respectively. If contact forms are 
selected as principle forms in the bundle of jets, then
$\omega_j^i =0$ and $\pi_j^i =0$ hold. The forms in \eqref{eqV11}
simplify to
\be
\tilde{\omega}_j^i = \Gamma^i_{jk} (x,y,u, \lambda_l) \, \omega^k,
\quad
\tilde{\pi}_j^i = \Phi^i_{jk} (x,y,v, \mu_l) \omega^k.
\label{eqV12}
\end{equation}
In this case, the equations $\theta=0$ and $\vartheta =0$
considered on pull-backs of solutions of the equations
\eqref{eqV6} and \eqref{eqV7}, respectively, are written as
\be
\begin{array}{c}
dv - \xi^i_j (v) \Gamma^j_{i1} (x,y,u, u_k) \, dx
- \xi^i_j (v) \Gamma^j_{i2} (x,y,u,u_k) \, dy =0,  \\
  \\
du - \eta_j^i (u) \Phi_{i1}^j (x,y,v,v_k) \, dx
- \eta^i_j (u) \Phi^j_{i2} (x,y,v,v_k) \, dy =0.  \\
\end{array}
\label{eqV13}
\end{equation}
Of course, \eqref{eqV13} are equivalent to the following
systems of partial differential equations
\be
v_x = \xi^i_j (v) \Gamma^j_{i1} (x,y,u,u_k),
\qquad
v_y = \xi^i_j (v) \Gamma^j_{i2} (x,y,u,u_k),
\label{eqV14}
\end{equation}
and 
\be
u_x = \eta^i_j (u) \Phi_{i1}^j (x,y,v, v_k),
\qquad
u_y = \eta^i_j (u) \Phi^j_{i2} (x,y,v,v_k).
\label{eqV15}
\end{equation}

\section{An Application of the Theory}

This formalism is now applied to obtain B\"acklund transformations
between the Liouville equation $u_{xy} =e^u$ and the wave equation
$v_{xy}=0$. These can now be defined by specifying the connections
in two principle bundles which define representations of zero
curvature, and the corresponding connections in the associated
bundles. In this case, the connection forms in the principle bundles
are defined as in \eqref{eqV12}.

A system of forms which will accomplish the task can be specified
as follows
\be
\begin{array}{cccc}
\tilde{\omega}^1_1 =- \dss\frac{\lambda_1}{4}  dx + \dss\frac{\lambda_2}{4}  dy, &
\tilde{\omega}^2_2 = \dss\frac{\lambda_1}{4} dx - \dss\frac{\lambda_2}{4} dy, &
\tilde{\omega}^2_1 = \dss\frac{1}{\sqrt{2}} e^{u/2}  dx, & 
\tilde{\omega}^1_2 = \dss\frac{1}{\sqrt{2}} e^{u/2} \, dy  \\
  &  &  &  \\
\tilde{\pi}^1_1 = \dss\frac{\mu_1}{4} dx - \dss\frac{\mu_2}{4} dy,  &
\tilde{\pi}_2^2 =- \dss\frac{\mu_1}{4} dx + \dss\frac{\mu_2}{4} dy, &
\tilde{\pi}^2_1 = \sqrt{2} (e^{-v/2} dx + e^{v/2} dy), & 
\tilde{\pi}^1_2 =0.  \\
\end{array}
\label{eqVI1}
\end{equation}
Based on this collection of definitions, the required coefficients
$\Gamma^j_{ik}$ and $\Phi^j_{ik}$ can be read off
\be
\begin{array}{cccc}
\Gamma^1_{11} =- \dss\frac{\lambda_1}{4}, & \Gamma_{21}^2 = \dss\frac{\lambda_1}{4}, &
\Gamma_{11}^2 = \dss\frac{1}{\sqrt{2}} e^{u/2}, & \Gamma^1_{21} =0,  \\
   &    &      &    \\
\Gamma^1_{12} = \dss\frac{\lambda_2}{4}, & \Gamma_{22}^2 =- \dss\frac{\lambda_2}{4}, &
\Gamma_{12}^2 =0,  & \Gamma^1_{22} = \dss\frac{1}{\sqrt{2}} e^{u/2}.  \\
\end{array}
\label{eqVI2}
\end{equation}
and as well,
\be
\begin{array}{cccc}
\Phi_{11}^1 = \dss\frac{\mu_1}{4}, & \Phi_{21}^2 =- \dss\frac{\mu_1}{4}, & 
\Phi_{11}^2 = \sqrt{2} e^{-v/2}, &  \Phi_{21}^1 =0,  \\
   &     &     &     \\
\Phi_{12}^1 =- \dss\frac{\mu_2}{4}, & \Phi_{22}^2 = \dss\frac{\mu_2}{4}, &
\Phi_{12}^2 = \sqrt{2} e^{v/2}, & \Phi_{22}^1 = 0.  \\
\end{array}
\label{eqVI3}
\end{equation}
Now the corresponding sets of forms $\theta^i$ and $\vartheta^i$ are
defined in terms of the structural forms \eqref{eqVI1},
\be
\begin{array}{ccc}
\theta^1 = 2 ( \tilde{\omega}^1_1 - \tilde{\omega}_2^2), &
\theta^2 = 2 \tilde{\omega}^2_1,  & \theta^3 =-2 \tilde{\omega}^1_2,  \\
    &     &     \\
\vartheta^1 = 2 ( \tilde{\pi}^1_1 - \tilde{\pi}^2_2 ),  &
\vartheta^2 = \tilde{\pi}_1^2,  &  \vartheta^3  =0.  \\
\end{array}
\label{eqVI4}
\end{equation}
It will be shown explicitly that these forms define representations of
zero curvature for the two equations above. The first three structure
equations in the $\theta^i$ are given by
$$
d \theta^1 + \theta^2 \wedge \theta^3 =- d \lambda_1 \wedge dx + d \lambda_2 \wedge dy
- 2 e^u \, dx \wedge dy,
$$
\be
d \theta^2 - \frac{1}{2} \theta^1 \wedge \theta^2 = \frac{1}{\sqrt{2}} e^{u/2}
(du - \lambda_2 dy) \wedge dx,
\label{eqVI55}
\end{equation}
$$
d \theta^3 + \frac{1}{2} \theta^1 \wedge \theta^3 = \frac{1}{\sqrt{2}} e^{u/2}
(-du + \lambda_1 \, dx) \wedge dy.
$$
On a section $\Sigma_1 \subset E_1$, using \eqref{eqII4} it follows that
$\lambda_1 =u_x$, $\lambda_2 = u_y$ and all three of these equations
vanish provided that $u$ satisfies $u_{xy}= e^u$. 

Similarly, for the forms $\vartheta^i$, it is found that
\be
\begin{array}{c}
d \vartheta^1 + \vartheta^2 \wedge \vartheta^3 = d \mu_1 \wedge dx - d \mu_2 \wedge dy,  \\
  \\
d \vartheta^2 - \frac{1}{2} \vartheta^1 \wedge \vartheta^2 =
\dss\frac{1}{\sqrt{2}} [- e^{v/2} (dv \wedge dx + \mu_2 dx \wedge dy)
+ e^{v/2} (dv \wedge dy - \mu_1 dx \wedge dy) ].
\end{array}
\label{eqVI66}
\end{equation}
The third vanishes identically since $\vartheta^3 =0$.
On a section $\Sigma_2 \subset E_2$, by applying \eqref{eqII4},
it follows that $\mu_1 = v_x$, $\mu_2 = v_y$ and these equations vanish
provided that $v$ satisfies the equation $v_{xy}=0$.

Using \eqref{eqVI2} and the definitions in \eqref{eqVI4}, then under the
assignment
\be
\xi_1^1 =2 , \quad \xi_2^2 =- 2,
\quad \xi_2^1 = 2 e^{-v/2},  \quad  \xi_1^2 =- 2 e^{v/2}.
\label{eqVI5}
\end{equation}
the equation for $dv$ in \eqref{eqV13} is written as
\be
dv - \theta^1 - e^{-v/2} \theta^2 - e^{v/2} \theta^3 =0.
\label{eqVI6}
\end{equation}
Similarly, using \eqref{eqVI3} and identifying
\be
\eta_1^1 =-2, \qquad \eta_2^2 =2, \qquad  \eta_2^1 = e^{u/2},
\label{eqVI7}
\end{equation}
the equation for $d u$ in \eqref{eqV13} becomes
\be
du + \vartheta^1 - e^{u/2} \vartheta^2 =0.
\label{eqVI8}
\end{equation}
It can be observed that the one-form of \eqref{eqVI6} is a
closed form, whereas \eqref{eqVI8} is not closed, but leads
to a consistent result.
Now all of the required information is at hand to
write down B\"acklund transformations between these two
equations. Substituting \eqref{eqVI2} and \eqref{eqVI5}
into \eqref{eqV14}, there results the system
\be
u_x + v_x = \sqrt{2} e^{(u-v)/2}, \qquad
u_y -v_y = \sqrt{2} e^{(u+v)/2}.
\label{eqVI9}
\end{equation}
Substituting \eqref{eqVI3} and \eqref{eqVI7} into
\eqref{eqV15}, it is found that the same pair appears,
\be
u_x + v_x = \sqrt{2} e^{(u-v)/2}, \qquad
u_y - v_y = \sqrt{2} e^{(u+v)/2}.
\label{eqVI10}
\end{equation}

{\bf Theorem 6.1.} The exterior derivatives of the one-forms
in \eqref{eqVI6} and \eqref{eqVI8} vanish modulo the sets
of forms $\{ dv, d \theta^i \}$ and $\{ du, d \vartheta^i \}$,
respectively.

{\bf Proof:} Let $\tau$ denote the one-form on the laft-hand side
of \eqref{eqVI6}. Differentiate $\tau$ exteriorly and there results,
$$
d \tau =- d \theta^1 + \frac{1}{2} e^{-v/2} dv \wedge \theta^2 -
e^{-v/2} d \theta^2 - \frac{1}{2} e^{v/2} dv \wedge \theta^3
- e^{v/2} d \theta^3.
$$
Replacing the known forms $d \theta^i$ from \eqref{eqVI55} and
$dv$ \eqref{eqVI6}, we obtain that
$$
d \tau = \theta^2 \wedge \theta^3 + \frac{1}{2} e^{-v/2} \theta^1 \wedge \theta^2
+ \frac{1}{2} \theta^3 \wedge \theta^2 - \frac{1}{2} e^{-v/2} \theta^1 \wedge
\theta^2 - \frac{1}{2} e^{v/2} \theta^1 \wedge \theta^3 - \frac{1}{2}
\theta^2 \wedge \theta^3 + \frac{1}{2} e^{v/2} \theta^1 \wedge \theta^3 =0,
$$
as required.

In the same way, differentiate \eqref{eqVI8} and substitute $d \vartheta^i$ from
\eqref{eqVI66} and $du$ from \eqref{eqVI8}.

This provides another way to get the $\xi^i_j$ and $\eta^i_j$
which appear in \eqref{eqV14} and \eqref{eqV15}.

{\bf Theorem 6.2.} Equations \eqref{eqVI9} and \eqref{eqVI10}
form a system of B\"acklund transformations which connect
the equations $u_{xy} = e^u$ and $v_{xy}=0$ respectively. 

{\bf Proof:} Differentiating the pair of equations in
\eqref{eqVI9} and replacing the first derivatives on the
right-hand side, it is found that
$$
(u + v)_{xy} = \frac{1}{\sqrt{2}} (u-v)_y e^{(u-v)/2} = e^u,
\qquad
(u- v)_{yx} = \frac{1}{\sqrt{2}} (u+v)_x e^{(u+v)/2} = e^u.
$$
Adding these two second derivatives, the Liouville 
equation $u_{xy} = e^u$ results.
Upon subtracting this pair, the wave equation $v_{xy} =0$ is obtained.

\section{Outlook and Summary}

A very general and useful formalism has been examined which
makes use of connections of zero curvature.
The first few sections present one way of giving an abstrct
formulation to this subject, and the latter part transfers this
to the more concrete aspect of actually calculating some
differential systems for a pair of specific equations.
If the forms are selected in the right way, it should be
possible to create auto-B\"acklund transformations, 
that is transformations between solutions of the same equation.
It has been shown that these types of connection have the
potential to produce Lax pairs and B\"acklund transformations
for nonlinear partial differential equations. In fact the
results of the previous section can be used to write Lax pairs
for the respective equations. Using coefficients 
\eqref{eqVI2} for the equation 
$u_{xy}=e^u$, the following Lax pair is obtained
$$
\begin{pmatrix}
X^1  \\
X^2  \\
\end{pmatrix}_x =
\begin{pmatrix}
\dss\frac{u_x}{4}  &  0  \\
- \dss\frac{1}{\sqrt{2}} e^{u/2} & - \dss\frac{u_x}{4} \\
\end{pmatrix} \begin{pmatrix}
X^1  \\
X^2  \\
\end{pmatrix},
\qquad
\begin{pmatrix}
X^1  \\
X^2  \\
\end{pmatrix}_y
= \begin{pmatrix}
- \dss\frac{u_y}{4}  & - \dss\frac{1}{\sqrt{2}} e^{u/2}  \\
0  &  \dss\frac{u_y}{4}  \\
\end{pmatrix} \begin{pmatrix}
X^1  \\
X^2  \\
\end{pmatrix}.
$$
The compatibility condition for this pair 
can be calculated by differentiating
the first matrix equation with respect to $y$ and the second with respect to $x$.
It is seen to hold provided that $u$ satisfies the equation $u_{xy}=e^u$.
Similarly, using the results in \eqref{eqVI3} for the
equation $v_{xy}=0$, the following Lax pair results
$$
\begin{pmatrix}
X^1  \\
X^2  \\
\end{pmatrix}_x =
\begin{pmatrix}
- \dss\frac{v_x}{4}  &  0  \\
- \dss\sqrt{2} e^{-v/2}  &  \dss\frac{v_x}{4}  \\
\end{pmatrix} \begin{pmatrix}
X^1  \\
X^2  \\
\end{pmatrix},
\qquad
\begin{pmatrix}
X^1  \\
X^2  \\
\end{pmatrix}_y = \begin{pmatrix}
\dss\frac{v_y}{4}  & 0  \\
- \dss\sqrt{2}  e^{-v/2}  & - \dss\frac{v_y}{4}  \\
\end{pmatrix}
\begin{pmatrix}
X^1  \\
X^2  \\
\end{pmatrix}.
$$
The compatibility condition is again found to hold
provided that $v$ satisfies $v_{xy}=0$.

It might be conjectured as a further application
of this work that if Lax pairs of the form
\eqref{eqV5} can be produced by some means, their matrix elements might
be used to generate 
connections of zero curvature as discussed here.
If they are found to have zero curvature structure, the results obtained
here would be of use in generating B\"acklund transformations
for the equations involved.

\newpage
\section{References}

\noindent
$[1]$ H. D. Wahlquist and F. B. Estabrook, Prolongation structures of nonlinear
evolution equations, J. Math. Phys. {\bf 16}, 1-7, (1975).  \\
$[2]$ F. B. Estabrook and H. D. Wahlquist, Prolongation structures of nonlinear
evolution equations II, J. Math. Phys. {\bf 17}, 1293-1297, (1976).  \\
$[3]$ F. B. Estabrook, Moving frames and prolongation algebras,
J. Math. Phys. {\bf 23}, 2071-2076, (1982).  \\
$[4]$ F. B. Estabrook, B\"acklund Transformations the Inverse Scattering Method, Solitons and 
Their Applications, Lecture Notes in Mathematics, ed. R. Miura, vol. 515, Springer, Berlin, 1976.  \\
$[5]$ R. Hermann, Pseudodifferentials of Estabrook and Wahlquist, the
geometry of solutions and the theory of connections,
Phys. Rev. Letts. {\bf 36}, 835-836, (1976).   \\
$[6]$ R. Hermann, The Geometry of NonLinear Differential Equations, B\"acklund
Transformations and Solitons, vol. XII, A, Math. Sci. Press, Brookline, MA, 1976.  \\
$[7]$ P. W. Michor, Topics in Differential Geometry, Graduate Studies in Mathematics,
vol. 93, AMS, Providence, RI, 2008.  \\
$[8]$ M. J. Ablowitz, D. K. Kaup, A. C. Newell and H. Segur, Nonlinear evolution
equations and physical significance, Phys. Rev. Letts.
{\bf 31}, 125-127, (1973).  \\
$[9]$ M. J. Ablowitz and H. Segur, Solitons and the Inverse Scattering Transform,
SIAM, Studies in Applied Mathematics, Philadelphia, PA, 1981.  \\
$[10]$ P. Bracken, A Geometric Interpretation of Prolongation by Means of
Connections, J. Math. Phys. {\bf 51}, 113502 (2010).  \\
$[11]$ P. Bracken, Exterior Differential Systems Prolongations and Application
to a Study of Two Nonlinear Partial Differential Equations,
Acta Appl. Math. {\bf 113}, 247-263, (2011).  \\
$[12]$ P. Bracken, Integrable Systems of Partial Differential Equations
Determined by Structure Equations and Lax Pair, Phys. Letts. {\bf A 374},
501-503, (2010).  \\
$[13]$ E. van Groesen and E. M. Jager, Mathematical structures in continuous
dynamical systems, Studies in Math. Physics, vol. 6, North Holland,
Amsterdam, II, Ch. 6, 1994.  \\
$[14]$ S. S. Chern and K. Tenenblat, Pseudospherical Surfaces and Evolution
Equations, Studies in Applied Math., {\bf 74}, 55-83, (1986).  \\
$[15]$ E. G. Reyes, Pseudo-spherical Surfaces and Integrability of
Evolution Equations, J. Diff. Equations, {\bf 147}, 195-230, (1998).  \\
$[16]$ I. M. Anderson and M. E. Fels, Symmetry Reduction of Exterior
Differential Systems and B\"acklund Transformations for PDE in the Plane,
Acta Appl. Math., {\bf 120}, 29-60, (2012).   \\
$[17]$ J. Krasilshchik and A. Verbovetsky, Geometry of Jet Spaces and
Integrable Systems, J. Geom. and Physics, {\bf 61}, 1633-1674, (2011).   \\
$[18]$ C. Rogers and W. K. Schief, B\"acklund and Darboux Transformations,
Cambridge Univ. Press, Cambridge, 2002.  \\
$[19]$ A. K. Rybnikov, Connections Defining Representations of Zero
Curvature and the Solitons of sine-Gordon and Korteweg-de Vries Equations,
Russian J. of Math. Phys. {\bf 18}, 195-210, (2011).  \\
$[20]$ A. K. Rybnikov, Equations of the Inverse Problem,
B\"acklund Transformations and the Theory of Connections,
J. of Math. Sciences, {\bf 94}, 1685-1699, (2009).  \\
$[21]$ R. L. Bryant, S. S. Chern, R. B. Gardner, H. L. Goldschmidt,
P. A. Griffiths, Exterior Differential Systems, Springer-Verlag, 1991.  \\
$[22]$ F. B. Estabrook and H. D. Wahlquist, Classical geometries defined by
exterior differential systems on higher frame bundles, 
Classical and Quantum Gravity {\bf 6}, 263-274, (1989).  \\
$[23]$ P. Bracken, Connections Defining Representations of Zero
Curvature and their Lax and B\"acklund Mappings, Journal of Geometry
and Physics, {\bf 70}, 157-163, (2013).  \\

\end{document}